
\input amstex
\documentstyle {amsppt}
\magnification = \magstep 1
\rm         
\hsize15.3truecm
\vsize22.5truecm
\NoBlackBoxes
\hoffset=5truemm
\voffset=-5truemm

 0

\def \Bl {\operatorname{Bl}}
\def \Ch {\operatorname{Ch}}
\def \Closure {\operatorname{Closure}}
\def \ord  {\operatorname{ord}}
\def \Sing {\operatorname{Sing}}
\def \Supp {\operatorname{Supp}}
\def \grad {\operatorname{grad \,}}
\def \dim {\operatorname{dim}}
\def \R {\operatorname{R}}

\def \bullz {\hskip2pt\strut^\centerdot\hskip-5pt z}
\def \hak {^{\vee}}
\def \hakP {\Bbb P\hskip-5.7pt\strut^{^\vee}}
\def \jeden {1\hskip-3.5pt1}

\def \kreska {\hbox{\vbox{\hrule width 4pt
                          \vskip2.3pt}}}
\def \przeryw {\,\kreska\hskip2pt\kreska\hskip2pt\kreska\hskip-6pt>}

\document
\topmatter
\title
 CHARACTERISTIC CLASSES OF HYPERSURFACES AND CHARACTERISTIC CYCLES
\endtitle
\rightheadtext{ADAM PARUSI\'NSKI AND PIOTR PRAGACZ}
\leftheadtext{CHARACTERISTIC CLASSES OF HYPERSURFACES}
\author 
Adam Parusi\'nski \ and \ Piotr Pragacz
\endauthor

\endtopmatter

\bigskip

{\eightrm
\subhead
Abstract
\endsubhead
We give a new formula for the Chern-Schwartz-MacPherson class of
a hypersurface, generalizing the main result of [P-P],
which was a formula for the Euler characteristic.
Two different approaches are presented.
The first is based on the theory of characteristic cycle and the work
of Sabbah [S], Brian\c con-Maisonobe-Merle [B-M-M],
and L\^e-Mebkhout [L-M].
In particular, this approach leads to a simple proof of a formula of Aluffi [A]
for the above mentioned class.
The second approach uses Verdier's [V] specialization property of the 
Chern-Schwartz-MacPherson classes.
Some related new formulas are also given.
}

\bigskip\bigskip

\head
{\bf Introduction and statement of the main result}
\endhead

\medskip

Let $X$ be a nonsingular compact complex analytic variety of pure
dimension $n$ and let $L$ be a holomorphic line bundle on $X$.
Take $f\in H^0(X,L)$ a holomorphic section of $L$ such that the variety $Z$
of zeros of $f$ is a (nowhere dense) hypersurface in $X$.
Recall, after [A], that the {\it Fulton class} of $Z$ is
$$
c^F(Z)=c(TX|_Z-L|_Z)\cap [Z]\,,
\tag 1
$$
where $TX$ denotes the tangent bundle of $X$.
Note that if $Z$ is nonsingular then $c^F(Z)=c(TZ)\cap [Z]$.
By $c_*(Z)$ we denote the {\it Chern-Schwartz-MacPherson class} of $Z$,
see [McP]. We recall its definition later in Section 1.
After [Y] we shall call
$$
\Cal M(Z)=(-1)^{n-1}\bigl(c^F(Z)-c_*(Z)\bigr)
\tag 2
$$
the {\it Milnor class} of $Z$. This class is supported on the singular locus
of $Z$; it is convenient, however, to treat it as an element of $H_*(Z)$.

\bigskip

\example{Example 0.1}
Suppose that the singular set of $Z$ is finite and equals
$\{x_1,\ldots,x_k\}$.
Let $\mu_{x_i}$ denote the Milnor number of $Z$ at $x_i$ (see [M]).
Then 
$$
\Cal M(Z)=\sum^k_{i=1}\mu_{x_i}[x_i]\in H_0(Z)
$$ 
- see, for instance Suwa [Su], where this result is generalized to complete 
intersections.
\endexample

\smallskip

Consider the function $\chi:Z\to\Bbb Z$ defined for $x\in Z$ by
$\chi(x):=\chi(F_x)$, where $F_x$ denotes the Milnor fibre at $x$
(see [M]) and $\chi(F_x)$ its Euler characteristic.
Define also the function $\mu:Z\to \Bbb Z$ by $\mu=(-1)^{n-1}(\chi-\jeden_Z)$.

Fix now any stratification $\Cal S=\{S\}$ of $Z$ such that 
$\mu$ is constant on the strata of $\Cal S$.  
For instance, any Whitney
stratification of $Z$ satisfies this property, see [B-M-M] or [Pa]. 
Actually, it is not difficult to see 
that the topological type of the Milnor fibres
is constant along the strata of Whitney stratification of $Z$.  
Let us denote the value of $\mu$ on the stratum $S$ by $\mu_S$.
Let
$$
\alpha(S)=\mu_S-\sum\limits_{S'\neq S,S\subset\overline{S'}}
\alpha(S')
\tag 3
$$
be the numbers defined inductively on descending dimension of $S$.
(These numbers appear as the coefficients in the development of $\mu$
as a combination of the $\jeden_{\overline S}$'s -- see Lemma 4.1.)

The main result of the present paper is

\proclaim{Theorem 0.2}
In the above notation,
$$
\Cal M(Z)=\sum\limits_{S\in\Cal S} \alpha(S) c(L|_Z)^{-1}\cap
(i_{\overline S,Z})_* c_*(\overline S),
\tag 4
$$
where $i_{\overline S,Z}:\overline S\to Z$ denotes the inclusion.
\endproclaim

\medskip

When $X$ is projective, (4) was conjectured in [Y].
Under this last assumption, the equality
$$
\int\limits_Z \Cal M(Z)=\sum\limits_{S\in\Cal S} \alpha(S)
\int\limits_{\overline S} c(L|_{\overline S})^{-1}\cap c_*(\overline S)
\tag 5
$$
was proved in [P-P]; hence the theorem gives, in particular, a generalization
of the main result (5) of [P-P] to compact varieties.

\medskip

Our proof of the theorem is based on a formula due to Sabbah [S],
which allows one to calculate the Chern-Schwartz-MacPherson class of
a subvariety in terms of the associated {\it characteristic cycle}.
In the case of hypersurface $Z$, this characteristic cycle was 
calculated in [B-M-M] and [L-M] in terms of the  
blow-up of the Jacobian ideal of a local equation of $Z$ in $X$.  
So the proof of Theorem 0.2 is obtained by putting this local 
description and the global data together, and expressing the characteristic
cycle of $Z$ in terms of the global blow-up of the singular subscheme
of $Z$.  
Here by the {\it singular subscheme} of $Z$ we mean 
the one defined locally by
the ideal $\left(f,{\partial f\over\partial z_1},\ldots,{\partial f \over
\partial z_n}\right)$, where $(z_1,\ldots,z_n)$ are local coordinates on $X$.

\smallskip

The approach used leads to a very simple proof of a formula for the 
Chern-Schwartz-MacPherson class of hypersurface in terms of some divisors
associated with the above blow-up.
This formula was originally obtained by Aluffi [A] by different methods.
Some new formulas for the Chern-Schwartz-MacPherson classes of the
constructible functions $\chi$ and $\mu$ are also given.

\smallskip

Finally, we show, using Verdier's specialization property of the 
Chern-Schwartz-MacPherson classes (see [V], and also [S] and [K2])
how to prove another conjecture of Yokura, which, combined with a result
from [Y], gives an alternative proof of Theorem 0.2 in the case of
projective $X$.
We find that this specialization argument somewhat better explains the essence
of the theorem.

\bigskip\bigskip

\head
{\bf 1. Chern-Mather classes and Chern-Schwartz-MacPherson classes}
\endhead

\medskip

We start by recalling some results of Sabbah [S].
Let for $X$ as in the introduction, $T\hak X$ denote the cotangent bundle
of $X$.
Let $V$ be an (irreducible) subvariety of $X$.
Denote by $c_M(V)$ (resp. $c_M\hak(V)$\,) the {\it Chern-Mather class}
of $V$ (resp. the {\it dual Chern-Mather class}).
Let us recall briefly their definitions.
Let $\nu:NB(V)\to V$ be the Nash blow-up of $V$.
By definition on $NB(V)$ there exists the ``Nash tangent bundle" $T_V$
which extends $\nu^*TV^0$, where $V^0$ is the regular part of $V$.
Define the following elements of $H_*(V)$
$$
\aligned
c_M(V) &= \nu_*\bigl(c(T_V)\cap[NB(V)]\bigr) \\
c_M\hak(V) &= \nu_*\bigl(c(T_V\hak)\cap[NB(V)]\bigr),
\endaligned
\tag 6
$$
where $T_V\hak$ is the dual bundle of $T_V$.
It is easy to see that
$$
c_M\hak(V)=(-1)^{\dim\,V} c_M(V)\hak ,
\tag 7
$$
where for a homology class $a=a_0+a_1+a_2+\ldots$ , where $a_i\in H_{2i}(V)$,
we denote $a\hak=a_0-a_1+a_2-\ldots$ .
\smallskip

By $T_V\hak X\subset T\hak X$ we denote the {\it conormal space} to $V$ :
$$
T_V\hak X= \Closure \left\{(x,\xi)\in T\hak X\,\mid\,x\in V^0, \ 
\xi|_{T_xV^0}\equiv 0\right\} ,
$$
and by $C(V)\subset \Bbb PT\hak X$ its projectivization.
Let $\pi:C(V)\to V$ be the restriction of the projection $\Bbb PT\hak X\to X$
to $C(V)$.
Then by [S], we have
$$
\aligned
c_M\hak(V) &= c(T\hak X|_V)\cap\pi_*\left(c\left(\Cal O(-1)\right)^{-1}\cap
[C(V)]\right) \\
c_M(V) &= (-1)^{n-1-\dim\,V} c\left(TX|_V\right)\cap\pi_*
\left(c\left(\Cal O(1)\right)^{-1}\cap [C(V)]\right) \,,
\endaligned
\tag 8
$$
where $\Cal O(-1)$ is the tautological line bundle on $\Bbb PT\hak X$
restricted to $C(V)$.

Let now $\varphi$ be a constructible function on $X$,
$$
\varphi=\sum a_j\jeden_{Y_j}\,,
$$
where $Y_j$ are (closed) subvarieties of $X$ and $a_j\in \Bbb Z$.
By the {\it characteristic cycle} of $\varphi$ we mean the Lagrangian
conical cycle in $T\hak X$ defined by
$$
\Ch(\varphi) = \Ch\left(\bigoplus\limits_j \left(i_{Y_j,X}\right)_*
\Bbb C\,\strut_{Y_j}^{\oplus a_j}\right)\,,
\tag 9
$$
where $\Bbb C\,\strut_{Y_j}$ is the constant sheaf on $Y_j$ and 
$i_{Y_j,X}:Y_j\to X$ denotes the inclusion.
For a general definition of the characteristic cycle of a sheaf, we refer
the reader to [B].
The characteristic cycle of a constructible function admits the following 
interpretation.
Let $F(X)$ and $L(X)$ denote the groups of constructible functions on $X$
and conical Lagrangian cycles in $T\hak X$ respectively.
It is known that the assignment
$$
T_V\hak X\mapsto (-1)^{\dim\,V} Eu_V\,,
\tag 10
$$
where $Eu_V$ stands for the Euler obstruction (see [McP] and also [S], [K1]),
defines a natural transformation of the functors of Lagrangian conical 
cycles and 
constructible functions, that is an isomorphism.
In particular, we have an isomorphism between $L(X)$ and $F(X)$.
The operation of taking the characteristic cycle is the inverse of this
isomorphism; that is, it is given by
$$
\Ch(Eu_V)=(-1)^{\dim\,V}T_V\hak X\,.
\tag 11
$$
Since every constructible function is a combination of the $Eu_V$'s
(see [McP]), this allows ``in theory" to compute $\Ch(\varphi)$
for a constructible function $\varphi$. However, even for $\varphi =
\jeden_V$, this would involve not only the Euler obstruction of $V$
itself but also of some subvarieties of $V$.

Now we associate with a constructible function $\varphi$ on $X$ its
{\it Chern-Schwartz-MacPherson class} (abreviation: CSM-class).
Let $\pi:\Supp\Bbb P\Ch(\varphi)\to \Supp \varphi$ be the restriction of the 
projection $\Bbb PT\hak X\to X$.
Set
$$
c_*(\varphi)=(-1)^{n-1}c\left(TX|_{\Supp \varphi}\right)\cap 
\pi_*\left(c\left(\Cal O(1)\right)^{-1}\cap[\Bbb P\Ch\varphi]\right)
\tag 12
$$
-- an element in $H_*(\Supp \varphi)$.
We note that, in particular, by (8), (11) and (12) one has
$$
c_*(Eu_V)=c_M(V)\,.
\tag 13
$$
If $V\subset X$ is a (closed) subvariety, we will write $c_*(V)=
c_*(\jeden_V)$ as is customary.
Note that (12) is in agreement with [McP] because for 
$\jeden_V=\sum_i b_i Eu_{Y_i}$, where $b_i \in \Bbb Z$ and 
$Y_i\subset X$ are (closed)
subvarieties, we have
$$
c_*(\jeden_V)=\sum\limits_ib_ic_*(Eu_{Y_i})=
\sum\limits_ib_ic_M(Y_i)=c_*(V)\,.
$$
Thus, denoting by $\pi:\Supp\Ch(\jeden_V)\to V$ the restriction of the
projection $\Bbb PT\hak X\to X$, we have
$$
c_*(V)=(-1)^{n-1}c\left(TX|_V\right)\cap
\pi_*\left(c\left(\Cal O(1)\right)^{-1}\cap[\Bbb P\Ch(\jeden_V)]\right)\,.
\tag 14
$$

\bigskip\bigskip

\head
{\bf 2. Characteristic cycle of a hypersurface (local case)}
\endhead

\medskip

Suppose that $U\subset\Bbb C\,\strut^n$ is an open subset and $Z\subset U$
is a hypersurface of zeros of a holomorphic function $f:U\to\Bbb C$\,.
Let $\Cal J_f$ denote the {\it Jacobian ideal} 
$\left({\partial f\over \partial z_1},\ldots,
{\partial f\over \partial z_n}\right)$ 
of $f$, where $(z_1,\ldots,z_n)$ are the standard coordinates of 
$\Bbb C\,\strut^n$.
Consider the blow-up $\pi:\Bl_{\Cal J_f}U\to U$ of $\Cal J_f$.
Recall that we may interpret it as follows
$$
\Bl_{\Cal J_f}U=\Closure\left\{ (x,\eta)\in U \times
\hakP\strut^{n-1}\,|\, x \notin \Sing Z,
\eta=\left[{\partial f\over\partial z_1}(x):\ldots:{\partial 
f\over \partial z_n}(x)
\right]\right\}\,,
$$
where $\Sing Z$ denotes the singular subscheme of $Z$.

\medskip

\remark{\bf Remark 2.1} \ $\Bl_{\Cal J_f}U$ can be also interpreted as the
projectivization of the {\it relative conormal space} 
$T_f\hak \subset T\hak U$ (see [B-M-M, \S\,2], where we put $\Omega=X=U$).
Then by the Lagrangian specialization all
fibres of the restriction of \ $\tilde f:T\hak U \to U @> f >> \Bbb C$ 
 \ to $T_f\hak$ are conical Lagrangian
subvarieties of $T\hak U$. In particular, every irreducible component of
$\tilde f^{-1}(0)\cap T_f\hak$ is conormal to its projection on $U$. 
For details we refer to [B-M-M, \S\,2] and to references therein.
\endremark

\bigskip

Let $\Cal Z$ be the total transform $\pi^{-1}(Z)$
of $Z$ in $\Bl_{\Cal J_f}U$ and $\Cal Z=\bigcup\limits_i D_i$ be the
decomposition of $\Cal Z$ into irreducible components.
Set $C_i=\pi(D_i)$ and denote by $\Cal I_{C_i}$ the ideal defining $C_i$.
Then define

\vbox{\vskip 10pt
\halign{\hskip 50pt $#$ \quad &=\quad multiplicity of \ $#$ along $D_i$ \cr
\hfill n_i &\hfil \Cal I_{C_i}\hfil \cr
\hfill m_i &\hfil f \hfil \cr
\hfill p_i &\hfil \Cal J_f \hfil \cr
}
\vskip 10pt}

Let us now record the following result.

\smallskip

\proclaim{Proposition 2.2} \hbox{\rm :}
\hskip 40pt $m_i=n_i+p_i$ .
\endproclaim

\demo{Proof} 
Observe that by Remark 2.1 we have $D_i=\Bbb PT_{C_i}\hak U$.
Let $x$ be a generic point of $C_i$ and choose a system of coordinates
$(z_1,\ldots,z_n)$ at $x$ such that $C_i=\{z_1=\ldots=z_k=0\}$
in a neighborhood of $x$.
Then, over a neighborhood of $x$, 
$$
D_i= C_i\times\hakP\strut^{k-1} ,
\tag 15
$$
where 
$$
\hakP\strut^{k-1}=\{[\eta_1:\ldots:\eta_n]\in\hakP\strut^{n-1}\,|\,
\eta_{k+1}=\ldots=\eta_n=0\}.
$$
Let $\zeta:E\to U$ denote the blow-up of the product of $\Cal J_f$ and
$\Cal I_{C_i}$\,. So
$$
E=\Closure \left\{\Bigl(x,[z_1(x):\ldots:z_k(x)],
\Bigl[{\partial f\over\partial z_1}(x):\ldots:{\partial 
f\over \partial z_n}(x)\Bigr]\Bigr)| x \notin \Sing Z\right \}
$$
in $U\times \Bbb P\strut^{k-1} \times \hakP\strut^{n-1}$.
Then $\zeta$ factors through $\pi$
$$
\aligned
E\hskip5pt \longrightarrow \hskip5pt  &  Bl_{\Cal J_f}U \\
\hskip 25pt
\vbox{\offinterlineskip
  \hbox{\hskip0.3pt$\diagdown$}
  \hbox{$\zeta\hskip4pt\diagdown$}
  \hbox{$\hskip17.2pt\searrow$} }
& \hskip 6.2pt
\vbox{\offinterlineskip
  \hbox{\hskip 1.3pt$|$}
  \hbox{\hskip 1.3pt$|$\hskip3pt$\pi$}
  \hbox{\hskip 0.4pt$\downarrow$}  }  \\
&\hskip 5pt U  
\endaligned
$$
and there exists at least one irreducible component, say $B_{ij}$, of the
exceptional divisor of $\zeta$ which projects surjectively onto $D_i$.
Let $\gamma(t)=\bigl(z(t),v(t),\eta(t)\bigr)$ be an analytic curve in $E$
such that $\bigl(z(0),v(0),\eta(0)\bigr)$ is a generic point of $B_{ij}$,
$z_{k+1}(t)\equiv\ldots\equiv z_n(t)\equiv 0$ and $f\bigl(z(t)\bigr)\ne 0$
for $t\ne 0$\,. Then we have for $t\ne 0$ 
$$
\aligned
v(t) & =[z_1(t):\,\ldots\,:z_k(t)]\in\Bbb P\strut^{k-1} \\
\eta(t) & =\left[{\partial f\over\partial z_1}\bigl(z(t)\bigr):\,\ldots\,:
{\partial f\over\partial z_n}\bigl(z(t)\bigr)\right]\in\hakP\strut^{n-1}
\endaligned
$$
and $\eta(0)=\bigl[\eta_1(0):\,\ldots\,:\eta_k(0):0:\,\ldots\,:0\bigr]$
by (15).

\medskip

Since $\bigl(z(0),\eta(0)\bigr)$ is a generic point of $D_{i}$ the 
following equality would imply the proposition 
$$
\aligned
&\ord_0(f\circ\zeta)\bigl(\gamma(t)\bigr)= \ord_0 f\bigl(z(t)\bigr)\\
&= \ord_0 \bigl(z_1(t),\ldots,z_k(t)\bigr)+
\ord_0 \left({\partial f\over\partial z_1}\bigl(z(t)\bigr),\ldots,
{\partial f\over\partial z_n}\bigl(z(t)\bigr)\right)\,.
\endaligned
\tag 16
$$

\smallskip

We show (16).  First we note that we
may suppose that $(z_1\circ\zeta,\ldots,z_k\circ\zeta)$ is generated by
$z_{i_0}\circ\zeta$ at $\gamma(0)$ and $\zeta^{-1}\Cal J_f$
is generated by ${\partial f\over \partial z_{j_0}}\circ\zeta$ at
$\gamma(0)$, where $j_0\in\{1,\ldots,k\}$) by (15).
We have
$$
\aligned
{d\over dt}f\bigl(z(t)\bigr) & =\sum\limits^k_{i=1}{\partial f\over
\partial z_i}\bigl(z(t)\bigr)\bullz_i(t)  \\
&= {\partial f\over\partial z_{j_0}}\bigl(z(t)\bigr)\cdot\bullz_{i_0}(t)
\left(\sum\limits^k_{i=1}{{\partial f\over\partial z_i}\bigl(z(t)\bigr)
\over {\partial f\over\partial z_{j_0}}\bigl(z(t)\bigr)} \cdot
{\bullz_i(t)\over\bullz_{i_0}(t)}\right)\,,
\endaligned
\tag 17
$$
where $\bullz_i$ stands for $\frac {dz_i}{dt}$.  
Note that the quotients make sense since $z_{i_0}\circ\zeta$ generates
$\zeta^{-1}(z_1,\ldots,z_k)$ and $\partial f/\partial z_{j_0}\circ\zeta$
generates $\zeta^{-1}\Cal J_f$.

\smallskip

We may suppose that $\eta_{j_0}=1$ and $v_{i_0}=1$, which corresponds
to choosing affine coordinates on $\Bbb
P\strut^{k-1}\times\hakP\strut^{n-1}$.  
Since 
$$
\lim_{t \to 0} \ [\bullz_1(t):\ldots : \bullz_k(t)] 
 = \lim_{t \to 0} \  [z_1(t): \ldots : z_k(t)]
$$  
we get 
$$
\lim_{t \to 0} \left(\sum\limits^k_{i=1}
{{\partial f\over\partial z_i}\bigl(z(t)\bigr)
\over {\partial f\over\partial z_{j_0}}\bigl(z(t)\bigr)} \cdot
{\bullz_i(t)\over\bullz_{i_0}(t)}\right)
= \lim_{t \to 0}
\left(\sum\limits^k_{i=1}{\eta_i(t) \over \eta_{j_0}(t)} \cdot
{v_i(t) \over v_{i_0}(t)} \right)
= \sum\limits^k_{i=1}\eta_i(0) v_i (0).  
$$
This last sum is nonzero by the transversality of relative polar
varieties, see, for instance,  [H-M, 8.7, Lemme de transversalit\'e].
Consequently, (17) implies 
$$
\ord_0 f\bigl(z(t)\bigr)-1=\ord_0 {\partial f\over\partial z_{j_0}}
\bigl(z(t)\bigr)+\bigl(\ord_0z_{i_0}(t)-1\bigr)
$$
which gives (16), as required.
\qed
\enddemo

\bigskip

In the following theorem, the equality (i) and the second equality in (ii) 
were established in [B-M-M] (see also [L-M]).

\bigskip

\proclaim{Theorem 2.3}
\vskip -25pt
$$
\aligned
\hskip 30pt \hbox{\rm (i) \ } 
& \Ch(\jeden_Z)=(-1)^{n-1}\sum_i n_i T_{C_i}\hak U \\
\hbox{\rm (ii) \ }
& \Ch(\chi)=\Ch(\R\Psi_f\Bbb C\,\strut_U)=(-1)^{n-1}\sum_i m_i T_{C_i}\hak U \\
\hbox{\rm (iii) \ } 
& \Ch(\mu)=(-1)^{n-1}\Ch\bigl(\R\Phi_f\Bbb C\,\strut_U\bigr)=
\sum_i p_i T_{C_i}\hak U
\endaligned
$$
\endproclaim
\noindent
For a definition of the complexes of nearby cycles $\R\Psi_f$ 
and vanishing cycles $\R\Phi_f$, we refer the reader
to [D-K].
The first equalities in (ii) and (iii) are well-known
and follow from the local index theorem, see for instance
[B-D-K] and [S, (1.3) and (4.4)].

\demo{Proof of (iii)}
By the definition of $\mu$ we have
$$
\Ch(\mu)=(-1)^{n-1}\bigl(\Ch(\chi)-\Ch(\jeden_Z)\bigr)\,.
$$
Hence, using Proposition 2.2, the assertion follows.
\qed
\enddemo

\bigskip

Let $\Cal Y$ denotes the exceptional divisor in $\Bl_{\Cal J_f}U$\,.
Since $D_i=\Bbb PT_{C_i}\hak U$, we can rewrite the assertions of the 
theorem as the following equalities.

\bigskip

\remark{\bf Corollary 2.4}
\vskip -25pt
$$
\aligned
\hskip 30pt \hbox{\rm (i) \ } 
& \left[\Bbb P\Ch(\jeden_Z)\right]=(-1)^{n-1}\left(\left[\Cal Z\right]-
\left[\Cal Y\right]\right) \\
\hbox{\rm (ii) \ }
& \left[\Bbb P\Ch(\chi)\right]=(-1)^{n-1}\left[\Cal Z\right] \\
\hbox{\rm (iii) \ } 
& \left[\Bbb P\Ch(\mu)\right]=[\Cal Y]
\endaligned
$$
\endremark

\smallskip

Observe that these equalities already take place on the level of cycles.

\smallskip

\remark{\bf Remark 2.5}
Since $f$ belongs to the integral closure of $\Cal J_f$
(see [LJ-T]) the normalizations of the blow-ups
of $\Cal J_f$ and $\left(f,{\partial f\over\partial z_1},\ldots,
{\partial f \over\partial z_n}\right)$ are equal.
Hence Corollary 2.4 holds true if we replace the blow-up of the former
ideal by the blow-up of the latter one.  
\endremark

\bigskip\medskip

\head
{\bf 3. Characteristic cycle of a hypersurface (global case)}
\endhead

\medskip

Let $X,\,L,\,Z,\,f$ be as in the introduction.
Let $B=\Bl_YX\to X$ be the blow-up of $X$ along the singular subscheme
$Y$ of $Z$.
Let $\Cal Z$ and $\Cal Y$ denote the total transform of $Z$ and the 
exceptional divisor in $B$, respectively.
The following description of the CSM-class of $Z$ was established by
Aluffi [A] by different methods.

\medskip

\proclaim{Theorem 3.1} {\rm ([A])
Let $\pi:\Cal Z\to Z$ be the restriction of the blow-up to $\Cal Z$.
Then
$$
c_*(Z)=c\left(TX|_Z\right)\cap
\pi_*\left({[\Cal Z]-[\Cal Y]\over 1+\Cal Z-\Cal Y}\right)\,,
$$
where on the RHS, $\Cal Z$ and $\Cal Y$ mean the first Chern classes of the
line bundles associated with $\Cal Z$ and $\Cal Y$ i.e. those of 
$\pi^*\left(L|_Z\right)$ and $\Cal O_B(-1)$, the latter being the canonical
line bundle on $B$.
\endproclaim

\demo{Proof}
To get a convenient description of $B$, we use (after [A]) the bundle
$\Cal P^1_XL$ of principal parts of $L$ over $X$ (see e.g. [At]).
Consider the section $X\to\Cal P^1_XL$ determined by $f\in H^0(X,L)$.
Recall that $\Cal P^1_XL$ fits in an exact sequence
$$
0\to T\hak X\otimes L\to\Cal P^1_XL\to L\to 0
$$
and the section in question is written locally as
$(df,f)=\left({\partial f\over\partial z_1},\ldots,{\partial f\over 
\partial z_n},f\right)$\,, where $(z_1,\ldots,z_n)$ are local coordinates
on $X$.
It follows that the closure of the image of the meromorphic map
$X\przeryw \Bbb P\Cal P^1_XL$ induced by $(df,f)$ is the blow-up
$B \to X$.
Thus we may treat $B$ as a subvariety of $\Bbb P\Cal P^1_XL$.
Clearly, the total transform $\Cal Z$ of $Z$ equals $B \cap\Bbb P(T\hak
X\otimes L)$.
The canonical line bundle $\Cal O_B(-1)=\Cal O(\Cal Y)$ on $B$ is the
restriction of the tautological line bundle $\Cal O(-1)$ on
$\Bbb P\Cal P^1_XL$.
Observe that the bundle $\Cal O(-1)$ restricted to $\Cal Z$
is contained in $(T\hak X\otimes L)|_{\Cal Z}$ (because $f\equiv 0$ over $Z$).
Hence $\Cal O_B(-1)|_{\Cal Z}$ is the restriction of the tautological
line bundle $\Cal O_{\widetilde{\Bbb P}}(-1)$ on 
$\widetilde{\Bbb P}=\Bbb P(T\hak X\otimes L)$.
Using the natural identification $\Bbb P(T\hak X\otimes L)\cong\Bbb P
(T\hak X)$ the line bundle $\Cal O_{\widetilde{\Bbb P}}(-1)$
corresponds to the line bundle $\Cal O_{\Bbb P}(-1)\otimes L$ on 
$\Bbb P=\Bbb P(T\hak X)$.
Thus $\Cal O_{\Bbb P}(1)$ on $\Bbb P$ corresponds to
$\Cal O_{\widetilde{\Bbb P}}(1)\otimes L$ on $\widetilde{\Bbb P}$.
Hence, using the characteristic cycle formula (14), we get
$$
\aligned
c_*(Z) &=(-1)^{n-1}c\left(TX|_Z\right)\cap\pi_*
\Bigl(c\bigl(\Cal O_B(1)\otimes \pi^*L|_Z\bigr)^{-1}\cap
\bigl[\Bbb P\Ch(\jeden_Z)\bigr]\Bigr) \\
& = c\left(TX|_Z\right)\cap\pi_*\left({[\Cal Z]-[\Cal Y]\over 1+\Cal Z-\Cal Y}
\right)
\endaligned
$$
because by (the global analogue of) Corollary 2.4, we have the equality
$[\Bbb P\Ch(\jeden_Z)]=(-1)^{n-1}\bigl([\Cal Z]-[\Cal Y]\bigr)$\,.\qquad
\qed
\enddemo

\bigskip

By Corollary 2.4, we have $[\Bbb P\Ch(\chi)]=(-1)^{n-1}[\Cal Z]$ and
$[\Bbb P\Ch(\mu)]=[\Cal Y]$.
Therefore, using similar arguments, we get the following result.

\smallskip

\proclaim{Theorem 3.2}
\vskip -25pt
$$
\aligned
\hskip 30pt \hbox{\rm (i) \ } 
& c_*(\chi)=c\left(TX|_Z\right)\cap
\pi_*\left({[\Cal Z]\over 1+\Cal Z-\Cal Y}\right)\,, \\
\hbox{\rm (ii) \ }
& c_*(\mu)=(-1)^{n-1}c\left(TX|_Z\right)\cap
\pi_*\left({[\Cal Y]\over 1+\Cal Z-\Cal Y}\right)\,.
\endaligned
$$
\endproclaim

\smallskip

\noindent
(The constructible function $\mu$ is supported on $Y$ 
but for later use we consider its CSM-class in $H_*(Z)$.)

\medskip

\remark{\bf Remark 3.3}
One can add to the above formulas also
$$
c_M(Z)=c_*(Eu_Z)=c\left(TX|_Z\right)\cap
\pi_*\left({[\Cal Z']\over 1+\Cal Z-\Cal Y}\right)\,,
$$
where $\Cal Z'$ is the proper transform of $Z$.
This equality for the Chern-Mather class was established originally by
Aluffi [A] by different methods.
Using the technique of characteristic cycles, it is a consequence of 
the equality
$\left[\Bbb P\bigl(\Ch(Eu_Z)\bigr)\right]=(-1)^{n-1}[\Cal Z']$ (see (11)).
\endremark

\bigskip\bigskip

\head
{\bf 4. Proof of Theorem 0.2}
\endhead

\smallskip

We start this section with the following fact about the constructible 
functions $\mu$ and $\alpha$ defined in the introduction.

\smallskip

\proclaim{Lemma 4.1} \hbox{\rm :}
\hskip80pt $\mu=\sum\limits_{S\in\Cal S}\alpha(S)\jeden_{\overline S}$\,.
\endproclaim

\demo{Proof}
Fix an arbitrary stratum $S_0$ and a point $x\in S_0$.
We have
$$
\aligned
&\biggl(\sum\limits_S \alpha(S)\jeden_{\overline S}\biggr)(x)=
\sum\limits_{S\neq S_0,\overline S\supset S_0} \alpha(S)+\alpha(S_0) \\
&= \sum\limits_{S\neq S_0,\overline S\supset S_0}\alpha(S) +
\left(\mu_{S_0}-\sum\limits_{S\neq S_0,\overline S\supset S_0}\alpha(S)\right)
=\mu(x)\,.\qquad \qed
\endaligned
$$
\enddemo

\medskip

Now we pass to the proof of Theorem 0.2.
Let $\pi:\Cal Z\to Z$ be the restriction of the blow-up $B=\Bl_YX\to X$.
We have, rewriting (1) as in [A] and using the projection formula,
$$
c^F(Z)=c\left(TX|_Z\right)\cap\pi_*\left({[\Cal Z]\over 1+\Cal Z}\right)\,.
$$
Invoking (2) and using Theorem 3.1, we get
$$
\aligned
\Cal M(Z) &= (-1)^{n-1}\bigl(c^F(Z)-c_*(Z)\bigr) \\
& = (-1)^{n-1}c\left(TX|_Z\right)\cap\pi_*\left({[\Cal Z]\over 1+\Cal Z}
-{[\Cal Z]-[\Cal Y]\over 1+\Cal Z-\Cal Y}\right) \\
& = (-1)^{n-1}c\left(TX|_Z\right)\cap\pi_*\left({[\Cal Y]\over (1+\Cal Z)
(1+\Cal Z-\Cal Y)}\right)
\endaligned
\tag 18
$$
because $\Cal Y\cap [\Cal Z]=\Cal Z\cap [\Cal Y]$.
If we pass to the characteristic cycle approach, the equality (18) is
rewritten, by Corollary 2.4, in the form
$$
\Cal M(Z)= (-1)^{n-1}c\left(TX|_Z\right)\cap\pi_*
\left({[\Bbb P\Ch(\mu)]\over (1+\Cal Z)(1+\Cal Z-\Cal Y)}\right)\,.
\tag 19
$$
Since $\mu=\sum_{S\in\Cal S}\alpha(S)\jeden_{\overline S}$ by Lemma 4.1,
we have 
$$
\Ch(\mu)=\sum_{S\in\Cal S}\alpha(S)\Ch(\jeden_{\overline S})
$$
and hence
$$
\aligned
&{[\Bbb P\Ch(\mu)]\over (1+\Cal Z)(1+\Cal Z-\Cal Y)}= \\
&=\sum\limits_{S\in\Cal S}\alpha(S)c(L|_Z)^{-1}\cap\pi_*
\Bigl(c\bigl(\pi^*L|_Z\otimes\Cal O_B(1)\bigr)^{-1}\cap[\Bbb P\Ch
(\jeden_{\overline S})]\Bigr).
\endaligned
\tag 20
$$
By (14) and the proof of Theorem 3.1, we get 
$$
\bigl(i_{\overline S,Z}\bigr)_* c_*(\overline S)=
(-1)^{n-1}c\bigl(TX|_Z\bigr)\cap\pi_*
\Bigl(c\bigl(\pi^*L|_Z\otimes\Cal O_B(1)\bigr)^{-1}\cap
[\Bbb P\Ch(\jeden_{\overline S})]\Bigr)
\tag 21
$$
for each stratum $S \in \Cal S$.
Finally, using (20) and (21), we rewrite (19) in the form
$$
\Cal M(Z)=\sum\limits_{S\in\Cal S} \alpha(S)\,c\bigl(L|_Z\bigr)^{-1}\cap
\bigl(i_{\overline S,Z}\bigr)_*\,c_*(\overline S)
$$
which is the required expression.
\qed

\bigskip\bigskip

\head
{\bf 5. Another approach via specialization}
\endhead

\medskip

In this section the setup is as in the introduction.
Additionally, we suppose that there exists
 a section $g\in H^0(X,L)$ such that $Z'=g^{-1}(0)$ 
is smooth and transverse to the strata of a (fixed) Whitney stratification
$\Cal S=\{S\}$ of $Z$.
For $t\in \Bbb C$ denote $f_t=f-tg$ and set $Z_t=f_t^{-1}(0)$.
In this section by $\Cal Z$ we will denote the following correspondence
in $X\times\Bbb C$\,:
$$
\Cal Z=\big\{(x,t)\in X\times\Bbb C \,|\,x\in Z_t\bigr\}\,.
$$
Denoting by $p:\Cal Z\to\Bbb C$ the restriction to $\Cal Z$ of the projection 
onto the second factor of $X\times\Bbb C$, we have $Z_t=p^{-1}(t)$ for
$t\in\Bbb C$.

Let $F(\Cal Z)$ (resp. $F(Z)$\,) denote the group of constructible
functions on $\Cal Z$ (resp. on $Z$).
Denote by
$$
\sigma_F:F(\Cal Z)\to F(Z_0=Z)
$$
the {\it specialization map of constructible functions} (see [V], [S] and
[K2], where a different notation is used).
Recall briefly its definition.
If $Y\subset\Cal Z$ is a (closed) subvariety, one sets for the generator
$\jeden_Y$,
$$
\bigl(\sigma_F\jeden_Y\bigr)(x)=\lim \limits_{t\to 0} \ 
\chi\bigl(B(x,\varepsilon)\cap Y_t\bigr)
$$
for any sufficiently small $\varepsilon>0$, where $B(x,\varepsilon)$
is the closed ball of radius $\varepsilon$ about $x$ and $Y_t=Y\cap Z_t$.
In our situation, we are aiming to compute $\sigma_F\jeden_{\Cal Z}$.
More explicitly, for $x\in Z$ we want to calculate
$$
(\sigma_F\jeden_{\Cal Z})(x)=\lim \limits_{t\to 0} \ 
\chi\bigl(B(x,\varepsilon)\cap Z_t\bigr)\,.
$$
This is the content of the following

\medskip

\proclaim{Proposition 5.1}
One has
$$
\bigl(\sigma_F\jeden_{\Cal Z}\bigr)(x)=\left\{
\chi(x)=1+(-1)^{n-1}\mu(x) \quad \hbox{for } x\not\in Z\cap Z'\hskip4pt
\atop
1 \hfill \hbox{ for } x\in Z\cap Z'\,.
\right.
$$
\endproclaim

\demo{Proof} 
If $x\not\in Z\cap Z'$ i.e. $g(x)\ne 0$, then
$$
Z_t=\bigl\{z\,|\, f(z)-tg(z)=0\bigr\}=\bigl\{z\,|\,f(z)/g(z)=t\bigr\}
$$
after restriction to a small ball is the Milnor fibre of $f/g$ at $x$,
and $f/g$ also defines $Z$ in a neighborhood of $x$.
The assertion follows.

\smallskip

Let now $x\in Z\cap Z'$.
We will use similar arguments to those used in Step 1 of the proof of 
Proposition 7 in [P-P].
Proceeding locally we can assume that $x$ is the origin in $\Bbb C\,\strut^n$,
that in our local coordinates $g(z)\equiv z_n$ and that $\{z_n=0\}$
is transverse to a fixed Whitney stratification $\Cal S=\{S\}$ of
$Z=\{f=0\}$.
Our goal is to show that for sufficiently small $\varepsilon>0$ and
$0<\delta<<\varepsilon$, if $t\in\Bbb C$ satisfies $0<|t|<\delta$, then
$$
Z_t\cap B_\varepsilon=\bigl\{(z_1,\ldots,z_n)\in\Bbb C\,\strut^n\,\Big|\,
|z|<\varepsilon,\,f-tz_n=0\bigr\}
$$
is contractible, where $B_\varepsilon=B(0,\varepsilon)$.
Set $V=\{f=z_n=0\}$.
If $\varepsilon$ is sufficiently small then $V\cap B_\varepsilon$ is
contractible.
So it suffices to retract $Z_t\cap B_\varepsilon$ onto
$V\cap B_\varepsilon$.
In what follows we shall proceed on $Z_t \smallsetminus V$ 
for $t$ sufficiently small.   
First note that since the stratification is Whitney and hence satisfies
the $a_f$ condition (see [B-M-M] or [Pa]), we have by the assumption
on 
transversality 
$$
\left|{\Bigl({\partial f\over\partial z_1},\ldots,{
\partial f\over\partial z_{n-1}}\Bigr)}
\right| \geq c\left|{\partial f\over \partial z_n}\right|
$$
for some universal $c>0$.
Therefore the linear forms $df(p)$ and $dz_n(p)$ are linearly
independent for $p\not\in\{f=0\}$.
So are clearly the forms $d(f-tz_n)$ and $dz_n$.  Consequently 
the orthogonal projection of  
$\grad |z_n|$ onto $Z_t=\{f-tz_n=0\}\smallsetminus V$ is nonzero, 
and we may normalize it so that the normalized vector field $\vec v$ 
satisfies
$$
\aligned
(i)\hskip15pt & {\partial|z_n|\over\partial\vec v}=1\ ; \\
(ii)\hskip15pt & {\partial(f-tz_n)\over\partial\vec v}=0 .
\endaligned
$$

\smallskip

We want, as well, the trajectories of this vector field do not leave 
$B_\varepsilon$.
For this we modify $\vec v$ near $S_\varepsilon=\bigl\{z\,\Big|\,
|z|=\varepsilon\}$.
Let $p\in V\cap S_\varepsilon$ and let $S$ be the stratum which contains $p$.
Let $p(s)\to p$ as $s\to 0$ be an analytic curve such that 
$f\bigl(p(s)\bigr)\ne 0$ for $s\ne 0$.
Then the limit $\eta$ of $df\bigl(p(s)\bigr)$ in $\hakP\strut^{n-1}$
as $s\to 0$, exists.
The forms $\eta$ and $dz_n$ are linearly independent by 
the assumption on transversality, and both vanish on the
tangent space to $S\cap\{z_n=0\}$.  Therefore, 
by the Whitney condition (b) for the closure of $S\cap\{z_n=0\}$, we get
the linear independence of $\eta$, $dz_n$ and $\sum_{i=1}^{n} z_idz_i$
at $p$.  Consequently, the orthogonal projection of  
$\grad |z_n|$ onto $S_\varepsilon \cap(Z_t\smallsetminus V)$ is nonzero 
in a neighborhood of $p$.  
Since $S_\varepsilon\cap V$ is compact, there exist a neighborhood $U$ of
$S_\varepsilon\cap V$ and a vector field $\vec w$ on
$U\smallsetminus\bigl(\{z_n=0\}\cup\{f=0\}\bigr)$ such that for $t$
small enough 
$$
\aligned
(i)\hskip15pt & {\partial|z_n|\over\partial\vec w}=1\ ; \\
(ii)\hskip15pt & {\partial(f-tz_n)\over\partial\vec w}=0 
 \ ; \\
(iii)\hskip15pt & {\partial\rho\over\partial\vec w}=0\hbox{ ,\qquad where }
\rho(z)=\parallel z\parallel^2\,.
\endaligned
$$
Using partition of unity we ``glue" $\vec w$ and $\vec v$ in order to get
a vector field $\vec u$ defined on $Z_t \smallsetminus V$ such that
$$
\aligned
(i)\hskip15pt & {\partial|z_n|\over\partial\vec u}=1\ ; \\
(ii)\hskip15pt & {\partial(f-tz_n)\over\partial\vec u}=0 \ ; \\
(iii)\hskip15pt & {\partial\rho\over\partial\vec u}=0
\hbox{\qquad on } S_\varepsilon.
\endaligned
$$
The flow of $\vec u$ allows us to retract $Z_t\cap B_\varepsilon$ onto
$Z_{t,c}=Z_t\cap B_\varepsilon\cap\{|z_n|\leq c\}$ for $c$ as small as we want.
On the other hand, for $c$ small enough, $Z_{t,c}$ retracts onto 
$V\cap B_\varepsilon=Z_t\cap B_\varepsilon\cap\{z_n=0\}$,
as required.
\qed
\enddemo

\smallskip

Now we want to pass to the {\it specialization map of homology classes}
$$
\sigma_H:H_*(Z_t)\to H_*(Z_0=Z)
$$
(see [V], [S] and [K2], where a different notation is used).
Recall briefly its definition.
Let $D\subset \Bbb C$ be a disk of a sufficiently small radius such that
the inclusion $Z=Z_0\subset p^{-1}(D)$ is a homotopy equivalence.
Thus for small nonzero $t\in D$ one defines the above  $\sigma_H$
as the composition
$$
H_*(Z_t) @> \ \ i_*\ \ >> H_*(p^{-1}D)\cong H_*(Z_0=Z)\,,
$$
where $i:Z_t\to p^{-1}D$ is the inclusion.
Recall now that Verdier's specialization property of CSM-classes asserts
the following.
For $\varphi\in F(\Cal Z)$ and $t$ sufficiently small, one has
$$
\sigma_Hc_*\bigl(\varphi|_{Z_t}\bigr)= c_*(\sigma_F\varphi)\,.
\tag 22
$$
(see [V] and also [S] and [K2]).

\medskip

Let us evaluate the both sides of (22) for $\varphi=\jeden_{\Cal Z}$\,.
The LHS reads simply $\sigma_Hc_*(Z_t)$.
As for the RHS, we have by Proposition 5.1
$$
\aligned
\sigma_F\jeden_{\Cal Z} &=\jeden_Z+(-1)^{n-1}\bigl(\mu\cdot
\jeden_{Z\smallsetminus Z\cap Z'}\bigr) \\
& = \jeden_Z+(-1)^{n-1}\bigl(\mu\cdot\jeden_Z-\mu\cdot\jeden_{Z\cap Z'}
\bigr)\,.
\endaligned
\tag 23
$$

Invoking the equality $\mu=\sum_S\alpha(S)\jeden_{\overline S}$
(see Lemma 4.1), Equation (23) is rewritten as
$$
\sigma_F\jeden_{\Cal Z}=\jeden_Z+(-1)^{n-1}
\left(\sum\limits_S\alpha(S)\jeden_{\overline S}-\sum\limits_S\alpha(S)
\jeden_{\overline S\cap Z'}\right)\,,
\tag 24
$$
and applying $c_*$ to (24) we get that the RHS of (22) is evaluated as
$$
\aligned
&c_*\bigl(\sigma_F\jeden_{\Cal Z}\bigr)= \\
&=c_*(Z)+(-1)^{n-1}\left\{\sum\limits_S\alpha(S)
\bigl[(i_{\overline S,Z})_*c_*(\overline S)-(i_{\overline S\cap Z',Z})_*
c_*(\overline S\cap Z')\bigr]\right\}\,,
\endaligned
$$
where $i_{\overline S\cap Z',Z}$ denotes the inclusion 
$\overline S\cap Z'\to Z$\,.

\smallskip

Suming up, by virtue of the specialization property (22), we have proved

\medskip

\proclaim{Proposition 5.2}
For the specialization map $\sigma_H:H_*(Z_t)\to H_*(Z)$, where $t\ne 0$
is small enough, one has
$$
\aligned
&\sigma_Hc_*(Z_t)= \\
&= c_*(Z)+(-1)^{n-1}\left\{\sum\limits_S\alpha(S)\bigl[(i_{\overline S,Z})_*
c_*(\overline S)-(i_{\overline S\cap Z',Z})_*c_*(\overline S\cap Z')\bigr]
\right\}\,.
\endaligned
$$
\endproclaim

\medskip

We now state the folowing result which appeared as Conjecture 1.12 in [Y].

\medskip

\proclaim{Theorem 5.3}
In the above notation, one has
$$
\Cal M(Z)=\sum\limits_{S\in\Cal S}\alpha(S)\bigl[(i_{\overline S,Z})_*
c_*(\overline S)-(i_{\overline S\cap Z',Z})_*c_*(\overline S\cap Z')\bigr]
$$
\endproclaim

\demo{Proof}
Observe that for $t$ like in Proposition 5.2, we have $c_*(Z_t)=c^F(Z_t)$
because $Z_t$ is smooth.
Moreover, since Fulton's class is expressed in terms of the Chern
classes of vector bundles, one has 
$\sigma_H\bigl(c^F(Z_t)\bigr)=c^F(Z)$.
We thus have
$$
\aligned
\Cal M(Z) &= (-1)^{n-1}\bigl(c^F(Z)-c_*(Z)\bigr) \\
&= (-1)^{n-1}\bigl(\sigma_Hc^F(Z_t)-c_*(Z)\bigr) \\
&= (-1)^{n-1}\bigl(\sigma_Hc_*(Z_t)-c_*(Z)\bigr) \\
&= \sum\limits_S\alpha(S)\bigl[(i_{\overline S,Z})_*c_*(\overline S)-
(i_{\overline S\cap Z',Z})_*c_*(\overline S\cap Z')\bigr]
\endaligned
$$
by Proposition 5.2.
\qed
\enddemo

\bigskip

Finally, arguing as in [Y,\S\,2] one shows that Theorem 5.3 implies
Theorem 0.2 when $X$ is projective.

\bigskip
\noindent
{\bf Acknowledgements.} 
We thank P. Aluffi and S. Yokura for sending us their preprints [A]
and [Y] in spring 1996 and summer 1997, respectively. 
These two preprints were very inspiring for us during the 
preparation of the present paper. Especially inspiring was
the letter of Aluffi (dated March 4, 1996) attached to [A] asking 
the second named author about the relationship between the approach 
taken in [A] and that of [P-P] in the context of 
Chern-Schwartz-MacPherson classes.
We hope that the present paper answers Aluffi's question.
\smallskip

During the preparation of this paper the second named author profited
the hospitality of the Technion in Haifa (Israel),  Universit\'e
d'Angers (France), and was partially supported by the 
KBN grant No. 2PO3A 02711.

\bigskip\bigskip

{\eightrm
\noindent
A.P. : D\'epartement de Math\'ematiques, Universit\'e d'Angers, 
2 Bd. Lavoisier, 49045 Angers Cedex 01, France

\noindent e-mail: parus\@tonton.univ-angers.fr

\bigskip
\noindent
P.P. : Mathematical Institute of Polish Academy of Sciences, Chopina 12,
87-100 Toru\'n, Poland

\noindent e-mail: pragacz\@mat.uni.torun.pl
}

\enddocument
\bye